\documentclass[12pt]{amsart}

\usepackage{amsfonts}

\newtheorem{theorem}{Theorem}[section]
\newtheorem{corollary}[theorem]{Corollary}
\newtheorem{lemma}[theorem]{Lemma}
\newtheorem{proposition}[theorem]{Proposition}
\theoremstyle{definition}

\newcommand{\cH}{\mathcal H}
\newcommand{\cBH}{\mathcal B(\mathcal H)}
\newcommand{\cKH}{\mathcal K(\mathcal H)}

\newcommand{\cM}{\mathcal M}
\newcommand{\cN}{\mathcal N}
\newcommand{\cL}{\mathcal L}
\newcommand{\e}{\varepsilon}
\newcommand{\codim}{{\rm codim}\,}

\def\NN{\mathbb{N}}

\begin{document}

\title[Product of two quasi-nilpotent operators]
{An operator is a product of two quasi-nilpotent operators if and only if it is not
semi-Fredholm}
\footnote{To appear in the Royal Society of Edinburgh Proceedings A (Mathematics).}

\author{Roman Drnov\v{s}ek}
\address{Department of Mathematics, University of Ljubljana, Jadranska 19, SI-1000
Ljubljana, Slovenia, e-mail : roman.drnovsek@fmf.uni-lj.si}

\author{Vladimir M\"{u}ller}
\address{Mathematical Institute, Czech Academy of Sciences, \v{Z}itn\'{a} 25, 115 67
Prague 1, Czech Republic, e-mail : muller@math.cas.cz }

\author{Nika Novak}
\address{Department of Mathematics, University of Ljubljana, Jadranska 19, SI-1000
Ljubljana, Slovenia, e-mail : nika.novak@fmf.uni-lj.si }


\begin{abstract}
We prove that a (bounded linear) operator acting on an infinite-dimensional,
separable, complex Hilbert space can be written as a product of two quasi-nilpotent
operators  if and only if it is not a semi-Fredholm operator. 
This solves the problem posed by  Fong and Sourour in 1984. 
We also consider some closely related questions. In particular, we show 
that an operator can be expressed as a product of two nilpotent operators 
if and only if its kernel and co-kernel are both infinite-dimensional.
This answers the question implicitly posed by Wu in 1989.

\end{abstract}

\maketitle

\section{Introduction}

Throughout the paper, let $\cH$ denote a complex, separable, infinite dimensional
Hilbert space.
Denote by $\cBH$ the algebra of all operators (i.e., bounded
linear transformations) on $\cH$ 
and by $\cKH$ the ideal of all compact operators. 
The essential spectrum $\sigma_e(T)$, left essential spectrum
$\sigma_{le}(T)$ and right 
essential spectrum $\sigma_{re} (T)$ of an operator $T \in \cBH$
is defined as the spectrum, left spectrum and right spectrum of
the class $T+\cKH$ in the Calkin algebra $\cBH/\cKH$.

Recall that an operator $T\in\cBH$ is called upper semi-Fredholm
if its range $T(\cH)$ is closed and $\dim\ker T<\infty$, and it is 
called lower semi-Fredholm if $\codim T(\cH)<\infty$ 
(then the range $T(\cH)$ is closed automatically).
One can show that $T$ is not upper semi-Fredholm if and only if for each
$\e>0$ and each subspace $\cM\subset\cH$ of a finite codimension
there exists a unit vector $x\in \cM$ such that $\|Tx\|<\e$.

It is well known that $\sigma_{e}(T)$, $\sigma_{le}(T)$ and
$\sigma_{re}(T)$ are the sets  of all complex numbers $\lambda$
such that $T-\lambda$ is not Fredholm, upper semi-Fredholm and
lower semi-Fredholm, respectively.
For the details of Fredholm theory 
we refer to the books \cite{Co90} and \cite{Mu03}.

In 1984, Fong and Sourour \cite{FS84} considered the question 
which operators on $\cH$ 
are products of two quasi-nilpotent operators. 
Their first observation was the following necessary condition for such operators.

\begin{proposition}\label{necessary}
If $T \in \cBH$ is a product of (two)  quasi-nilpotent operators, then 
$0 \in \sigma_{e}(T^* T) \cap \sigma_{e} (T T^*)$,
or equivalently $0 \in \sigma_{le}(T) \cap \sigma_{re} (T)$. 
\end{proposition}

In other words, a necessary condition for such operators is that $T$ is not a
semi-Fredholm operator.
Proposition \ref{necessary} follows from the following assertion
that is an easy
consequence 
of the fact that a quasi-nilpotent element in a unital Banach algebra is 
a topological divisor of zero (see e.g. \cite[section XXIX.3]{GG93}). 

\begin{proposition} \label{divisor}
Let $t$ and $q$ be elements in a unital Banach algebra. If $q$ is quasi-nilpotent,
then 
$t q$ is not left invertible, and $q t$ is not right invertible. 
\end{proposition}

Fong and Sourour found the following sufficient condition for an operator to be the
product of two quasi-nilpotent operators (see \cite[Theorems 5 and 6]{FS84}).

\begin{theorem}\label{sufficient}
Let $T \in \mathcal B(\mathcal H)$. 
If $0 \in \sigma_{e} (T^* T + T T^*)$, then $T$
is a product of 
two quasi-nilpotent operators. Moreover, if $T$ is a
compact operator, then 
$T$ is a product of two compact quasi-nilpotent operators.
\end{theorem}

The authors of  \cite{FS84} left open the question whether the necessary condition 
from Proposition \ref{necessary} is also sufficient for an operator to be the
product of two quasi-nilpotent operators. 
The main result of our paper, proved in Section 2, gives an affirmative
answer to this question. In Section 3 we consider a question of finding 
common quasi-nilpotent factors for a given finite collection of operators.
Finally, in Section 4 we study similar problems replacing quasi-nilpotent operators
by nilpotent ones. We show that an operator $T$ can be expressed as a product 
of two nilpotent operators if and only if $\dim \ker T=\dim \ker T^* = \infty$.
This answers the question implicitly posed in the survey paper \cite[p.55]{Wu89}.
\\

\section{Products of two quasi-nilpotent operators}

We start with the following description of a non-semi-Fredholm operator. 

\begin{proposition}\label{p1}
Suppose that $T \in \cBH$ is not a semi-Fredholm operator. Then $\mathcal H$ can be
decomposed 
as a direct sum of three infinite-dimensional closed subspaces,  so that $T$ 
is similar to an operator of the form 
 $$\left(\begin{matrix}
 0 & A &0\\
 K & C &D\\
 0 & L &0
 \end{matrix} \right),$$
where $K$ and $L$ are compact operators.
\end{proposition}

\begin{proof}

Since neither $T$ nor $T^*$ is upper semi-Fredholm, we can find
inductively an orthonormal sequence $f_1,g_1,f_2,g_2,\dots$ in
$\cH$ such that
$\|Tf_n\|\le\frac{1}{n}$ and $\|T^*g_n\|\le\frac{1}{n}$.

Let $\mathcal M$ be the closed linear span of the vectors
$\{f_n\}$,  and let $\mathcal N$
be the closed 
linear span of $\{g_n\}$. 
Then ${\mathcal H} = {\mathcal M} \oplus {\mathcal L} \oplus
{\mathcal N}$, where 
${\mathcal L} = \left({\mathcal M} \oplus {\mathcal N} \right)^{\bot}$. 
Reducing ${\mathcal M}$ and ${\mathcal N}$ (if necessary) we may assume that the 
subspace ${\mathcal L}$ is infinite-dimensional.
The matrix of $T$ relative to this decomposition is
$$\left(\begin{matrix}
K_1&A&B \\ K_2&C&D \\ K_3&K_4&K_5
\end{matrix}\right),$$
where $K_1, \ldots, K_5$ are Hilbert-Schmidt operators, and hence compact.

We now claim that the subspaces $\cM$ and $\cN$ in the decomposition 
${\mathcal H} = {\mathcal M} \oplus {\mathcal L} \oplus  {\mathcal N}$ can be
reduced to  
suitable infinite-dimensional closed subspaces (and so the subspace ${\mathcal L}$
is enlarged)
so that the corresponding operator $B$ is equal to $0$. 
To end this, we choose decompositions $\cN = \cN_1 \oplus \cN_2$ and $\cM = \cM_1
\oplus \cM_2$
such that $\cN_2$ and $\cM_1$ are infinite-dimensional closed subspaces and 
$B(\cN_2) \subseteq \cM_2$. Then $B$ as an operator from $\cN = \cN_1 \oplus \cN_2$ 
to $\cM = \cM_1 \oplus \cM_2$ has the form 
$\left(\begin{matrix}
*&0\\ *&*
\end{matrix}\right).
$
Hence the operator $T$ relative to the decomposition 
$\cH = \cM_1 \oplus (\cM_2 \oplus \cL \oplus \cN_1) \oplus \cN_2$ 
has the upper right corner equal to $0$. 
This shows what we have claimed. So, we may assume that $B=0$. 

Using the same argument we can also show that there is no loss of generality in
assuming 
that $K_3=0$. 

In order to complete the proof, we apply \cite[Theorem 2]{AS71} to
show that 
the operators $K_1$ and $K_5$ are similar to operators of the form 
$\left(\begin{matrix}
0&* \\ *&*
\end{matrix}\right)$ 
and 
$\left(\begin{matrix}
*&* \\ *&0
\end{matrix}\right)$
respectively, where both zeroes act on infinite-dimensional Hilbert spaces.
This fact easily implies the desired conclusion.
\end{proof}

As in \cite{FS84}, the following lemma is useful for our study.

\begin{lemma}\label{L9}
If $A$ and $C$ are quasi-nilpotent operators and $B$ is any operator, then 
$\left(\begin{matrix}
A&B \\0&C
\end{matrix}\right)$ is quasi-nilpotent.
\end{lemma}

Now we state the main result of this paper.
\begin{theorem} \label{two}
An operator $T \in \cBH$ is a product of two quasi-nilpotent operators if and only if
it is not a semi-Fredholm operator. 
\end{theorem}

\begin{proof}
By Proposition \ref{necessary} the condition is necessary. To
show that it is also
sufficient,
assume that $T$ is not a semi-Fredholm operator.
By Proposition \ref{p1}, we may assume that $T$ can be represented as 
$$\left(\begin{matrix}
0&A&0\\ K&C&D\\ 0&L&0
\end{matrix}\right) $$
on $\cH = \cH_1 \oplus \cH_2 \oplus \cH_3$, where each $\cH_j$ is isomorphic to
$\cH$, and 
$K$ and $L$ are compact operators.

Since $\cH_1$ and $\cH_3$ are infinite-dimensional, we can write each of them as an
infinite sum of Hilbert spaces that are isomorphic to $\cH$ so that in the obtained
decomposition of the space $\cH$ the operator $T$ can be represented as 
$$
{\footnotesize  \left(\begin{matrix}
      &   &\vdots&\vdots  &\vdots&   &\\
      &   &0     &A_2     &0     &   &\\
\cdots&0  &0     &A_1     &0     &0  &\cdots\\
\cdots&K_2&K_1   &\fbox{C}&D_1   &D_2&\cdots\\
\cdots&0  &0     &L_1     &0     &0  &\cdots\\
      &   &0     &L_2     &0     &   &\\
      &   &\vdots&\vdots  &\vdots&   &
\end{matrix}\right)},
$$
where $\{K_n\}$ and $\{L_n\}$ are compact operators satisfying 
$\max\{ \Vert K_n\Vert, \Vert L_n \Vert\} \le 4^{-n}$ for all $n \ge 2$. 
Now define the operators $Q_1$ and $Q_2$ on $\cH$ by 
$$Q_1= 
{\footnotesize  \left(\begin{matrix}
      &       &       &      &\vdots  &\vdots&\vdots&      &\\
      &       &       &      &0       &A_3   &0     &      &\\
      &       &       &      &0       &A_2   &0     &      &\\
      &       &       &      &0       &A_1   &0     &0     &\cdots\\
\cdots&2^{-2}I&2^{-1}I&I     &\fbox{0}&C     &0     &0     &\cdots\\
      &       &       &      &0       &0     &2L_1  &      &\\
      &       &       &      &0       &      &0     &2^2L_2&\\
      &       &       &      &0       &      &      &0     &\ddots\\     
      &       &       &      &\vdots  &      &      &      &\ddots      
\end{matrix}\right)}
$$
and
$$Q_2=
{\footnotesize
\left(\begin{matrix}
\ddots&\ddots &       &      &\vdots  &\vdots&   &      &\\
      &0      &2^2K_2 &      &0       &0     &   &      &\\
      &       &0      &2K_1  &0       &0     &0  &0     &\cdots\\
      &       &       &0     &0       &D_1   &D_2&D_3   &\cdots\\
\cdots&0      &0      &0     &\fbox{0}&0     &0  &0     &\cdots\\
      &       &       &      &I       &      &   &      &\\
      &       &       &      &2^{-1}I &      &   &      &\\
      &       &       &      &2^{-2}I &      &   &      &\\     
      &       &       &      &\vdots  &      &   &      &      
\end{matrix}\right)}.
$$
To prove that the operator $Q_1$ (and similarly, $Q_2$) is quasi-nilpotent, it
suffices, by Lemma \ref{L9}, to show that the operator 
$$ R = 
{\footnotesize  \left(\begin{matrix}
 0     &  2 L_1    &0       & \cdots \\
 0     &   0       & 2^2L_2 & \cdots \\
 0     &   0       &    0   & \cdots \\
\vdots & \vdots    & \vdots & \ddots            
\end{matrix}\right)}
$$
is quasi-nilpotent. But this follows easily from the following estimate
\begin{eqnarray*}
\Vert R^n \Vert&=& \sup \{ \Vert (2^{j+1} L_{j+1})(2^{j+2} L_{j+2}) \ldots
(2^{j+n} L_{j+n}) \Vert : j = 0, 1, 2, \ldots\} \cr
&\le& \max\{\Vert 2^j L_j\Vert: j\in\NN\} \cdot 2^{-(2 + 3 + \ldots + n)} \ . 
\end{eqnarray*}
Since a direct computation shows that $T = Q_1 Q_2$, the proof is complete.  
\end{proof}

Theorem \ref{two} has an interesting connection with the theory of integral operators. 

\begin{corollary} \label{integral}
Let $T$ be an operator on $L^2(X, \mu)$, where $(X, \mu)$ is a
$\sigma$-finite measure space which is not completely atomic. Then the following
assertions 
are equivalent: \\
(a) The operators $T$ and $T^*$ are unitarily equivalent to integral operators on
$L^2(X, \mu)$. \\
(b) The operator $T$ is a product of two quasi-nilpotent operators. \\
(c) The operator $T$ is not a semi-Fredholm operator.
\end{corollary}

\begin{proof}
Recall \cite[section 15]{HS78} that an operator $A$ on $L^2(X, \mu)$ is unitarily
equivalent to an integral operator on $L^2(X, \mu)$ if and only if $0 \in
\sigma_{re}(A)$. 
Therefore, the assertion (a) holds if and only if 
$0 \in \sigma_{le}(T) \cap \sigma_{re} (T)$, i.e.,  
$T$ is not a semi-Fredholm operator. Thus, the assertions (a) and (c) are equivalent.
The equivalence of (b) and (c) was shown in Theorem \ref{two}.
\end{proof}

If in Corollary \ref{integral} the operators $T$ and $T^*$ are simultaneously
unitarily 
equivalent to integral operators on $L^2(X, \mu)$, i.e., there exists a unitary
operator $U$ 
on $L^2(X, \mu)$ such that both $UTU^*$ and $UT^*U^*$ are integral operators, 
then we have $0\in \sigma_e(TT^*+T^*T)$, so that even the assumption of Theorem
\ref{sufficient}
is satisfied. However, the unitary equivalence in Corollary
\ref{integral}
is not necessarily simultaneous. Namely, it is shown in \cite[Example 15.10]{HS78} that
there exists an operator $A$ such that both $A$ and $A^*$ are unitarily equivalent to
integral
operators, but $A + A^*$ is invertible. It follows that  the operators $A$ and $A^*$
cannot be
simultaneously unitarily equivalent to integral operators on $L^2(X, \mu)$.  \\

\section{Common quasi-nilpotent factors}

In this section we consider the following related problem. 
Given operators $T_1, \ldots, T_n$ on $\cH$,
we are searching for quasi-nilpotent operators $Q_1$ and $Q_2$ and 
operators $S_1, \ldots, S_n$ such that $T_i=Q_1 S_i Q_2$ for all $i=1, \ldots,n$. 
Using Proposition \ref{divisor} it is easy to verify that a necessary
condition for these factorizations is that  
$0 \in \sigma_e(\sum_{i=1}^{n} T_iT^*_i) \cap \sigma_e(\sum_{i=1}^{n} T^*_iT_i)$.
Theorem \ref{c10} below shows that this condition is also sufficient. 

We first consider the case when $T_1, \ldots, T_n$ are compact operators. 

\begin{proposition}\label{p4}
 For compact operators $K_1, \ldots, K_n$ on $\cH$ the following assertions hold.
 \begin{enumerate}
 \item[(a)] There exist a quasi-nilpotent operator $Q$ and compact operators $L_1,
\ldots, L_n$ such that
 $$K_i = L_i Q \, ,\quad i=1, \ldots, n.$$
 \item[(b)] There exist a quasi-nilpotent operator $Q$ and compact operators $L_1,
\ldots, L_n$ such that
 $$K_i = Q L_i  \, ,\quad i=1, \ldots, n.$$
 \item[(c)] There exist quasi-nilpotent operators $Q_1$ and $Q_2$ and compact operators
$L_1, \ldots, L_n$ such that
 $$K_i = Q_1 L_i Q_2  \, ,\quad i=1, \ldots, n.$$
 \end{enumerate}
\end{proposition}

\begin{proof}
First we prove part (a).
Since $K = \sum_{i=1}^n K_i^* K_i$ is a positive compact operator, there exists  an
orthonormal basis
$\{\varphi_j\}_{j \in\NN}$ with respect to which $K$ has a
diagonal form ${\rm
diag} \, (\lambda_1, \lambda_2, \ldots)$,
where $\{\lambda_j\}_{j \in \NN}$ is a decreasing sequence of non-negative numbers
converging to $0$.

Define the operator $Q$ on $\cH$ by 
$$Q \varphi_j = \sqrt[4]{\lambda_j} \, \varphi_{j+1} \text{ for all } j \in \NN.$$
To show that $Q$ is quasi-nilpotent, we simply notice that $\Vert Q^{2n} \Vert =
\sqrt[4]{\lambda_1 \lambda_2 \ldots
\lambda_{2n}} \le (\lambda_1 \lambda_{n+1})^{\frac{n}{4}}$, and therefore
$\displaystyle{\lim_{n\to \infty}
\Vert Q^n\Vert^{\frac{1}{n}} =0}$.

Now define linear transformations $L_1,\ldots, L_n$ on $\cH$ by 
$$L_i \varphi_1=0 \quad\text{ and }\quad L_i \varphi_j =
\frac{1}{\sqrt[4]{\lambda_{j-1}}}K_i \varphi_{j-1}
 \text{ for } j >1.$$
Since we have, for all $j > 1$,
$$\Vert L_i \varphi_j \Vert^2 = 
\frac{1}{\sqrt{\lambda_{j-1}}} \langle K_i^* K_i \varphi_{j-1}, \varphi_{j-1}
\rangle \le
\frac{1}{\sqrt{\lambda_{j-1}}} \langle K \varphi_{j-1}, \varphi_{j-1} \rangle 
\le \sqrt{\lambda_{j-1}},$$
we conclude that $L_1,\ldots, L_n$ belong to $\cBH$.
Define the operator $L$ on $\cH$ by  
$$L = \sum_{i=1}^n L_i^*L_i =  {\rm diag} \, (0, \sqrt{\lambda_1}, \sqrt{\lambda_2},
\ldots).$$
Since $\displaystyle{\lim_{j\to \infty}  \lambda_j =0}$, the operator $L$ is compact.

Since the ideal $\cKH$ is hereditary in $\cBH$ (see, e.g., \cite[Theorem I.5.3]{Da96}),
the operators $L_i^*L_i$ are compact, and so are the square roots
$(L_i^*L_i)^{1/2}$. Using the polar
decomposition we conclude that the operators $L_1, \ldots, L_n$ are
compact as well.
Since $L_i Q = K_i$ for all $i$, the assertion (a) is proved.

For the proof of part (b) we apply part (a) for the operators  $K_1^*, \ldots, K_n^*$ 
and then take the adjoints.
Since part (c) is an immediate consequence of parts (a) and (b), the proof is
complete. 
\end{proof}

Part (c) in the following result can be considered as 
a generalization of \cite[proposition 4]{FS84}.

\begin{theorem}\label{c10}
For operators $T_1, \ldots, T_n$ on $\cH$ the following assertions hold.
\begin{enumerate}
\item[(a)] If $0 \in \sigma_e(\sum_{i=1}^{n} T_iT^*_i)$, then there exist
a quasi-nilpotent operator
$Q$ and operators $S_1, \ldots, S_n$ on $\cH$ such that 
$$T_i =Q S_i, \quad i=1, \ldots, n.$$
\item[(b)] If $0 \in \sigma_e(\sum_{i=1}^{n} T^*_iT_i)$, then there exist
a quasi-nilpotent operator
$Q$ and operators $S_1, \ldots, S_n$ on $\cH$ such that 
$$T_i =S_i Q, \quad i=1, \ldots, n.$$
\item[(c)] If $0 \in \sigma_e(\sum_{i=1}^{n} T_iT^*_i) \cap \sigma_e(\sum_{i=1}^{n}
T^*_iT_i)$, 
then there exist quasi-nilpotent operators $Q_1$ and $Q_2$ and
operators $S_1, \ldots, S_n$ with $0 \in \sigma_e(\sum_{i=1}^{n} S_iS^*_i) \cap
\sigma_e(\sum_{i=1}^{n} S^*_iS_i)$ such that 
$$T_i =Q_1 S_i Q_2, \quad i=1, \ldots, n.$$
\end{enumerate}
\end{theorem}

\begin{proof}
Let $T=\sqrt{\sum_{i=1}^{n} T_iT^*_i}$. Since $0 \in \sigma_e(T)$, by
\cite[Proposition 2]{FS84} there exists a quasi-nilpotent operator $Q$ on $\cH$ 
such that $T^2 = Q Q^*$. Since $T_iT^*_i \le T^2 = Q Q^*$ for $i=1, \ldots, n$, 
it follows from the well-known theorem of Douglas \cite{D66} that there
exist operators $S_1, \ldots, S_n$ such that $T_i = QS_i$ for $i=1, \ldots, n$.
This completes the proof of part (a). Part (b) follows from part (a) by duality. 

Let us now prove part (c). 

We claim that there 
exist an invertible operator $V$ on $\cH$ and a decomposition of
the space $\cH$ such that 
$$ 
VT_i V^{-1}= 
\left(\begin{matrix}
0&A_i&0\\K_i&C_i&D_i\\0&L_i&0
\end{matrix}\right),
$$
where $K_i$ and $L_i$ are compact operators for $i=1, \ldots, n$.

Let $A=\sqrt{\sum_{i=1}^{n} T^*_i T_i}$ and
$B=\sqrt{\sum_{i=1}^{n} T_i T^*_i}$. 
Since the selfadjoint operators $A$ and $B$ are not Fredholm, they 
are not upper semi-Fredholm.
Therefore we can construct inductively an orthonormal sequence
$f_1,g_1,f_2,g_2,\dots$ in $\cH$ such that
$\|Af_k\|\le\frac{1}{k}$ and $\|Bg_k\|\le\frac{1}{k}$ for all
$k$. 
Note that for $i=1,\dots,n$ and $k\in\NN$ we have
$\|T_if_k\|\le\|Af_k\|\le\frac{1}{k}$ and
$\|T^*_if_k\|\le\|Bf_k\|\le\frac{1}{k}$.

Let $\cM$ be the closed span of the vectors $f_1,f_2,\dots$, and
$\cN$ the closed span of the vectors $g_1,g_2,\dots$. Let $\cL$
be the orthogonal complement of $\cM\oplus\cN$.
Then in the decomposition 
$\cH= \cM \oplus \cL \oplus \cN$ the operators $\{T_i\}_{i=1}^n$ 
are of the form 
$$T_i=\left(\begin{matrix}
K_i^{(1)}&A_i&B_i\\K_i^{(2)}&C_i&D_i\\K_i^{(3)}&K_i^{(4)}&K_i^{(5)}
\end{matrix}\right),
$$
where $K_i^{(j)}$ are compact operators. 

We now claim that the decomposition of the space can be chosen in such a way 
that all $B_i$ are equal to $0$. As in the proof of Proposition \ref{p1} we may
assume that $B_1=0$. 
The subspaces $\cM$ and $\cN$ can be reduced to infinite-dimensional closed subspaces 
$\cM_1$ and $\cN_1$ so that in the decomposition of the space 
$\cH= \cM_1 \oplus \cL_1 \oplus \cN_1$ both operators $B_1$ and $B_2$ are equal to
$0$. 
After finitely many reductions of $\cM$ and $\cN$ we obtain the desired
decomposition.

Using the same argument as in Proposition \ref{p1} we can now find  
an invertible operator $V$ on $\cH$ and a decomposition
of the space such that
$$V T_i V^{-1} =\left(\begin{matrix}
0&A_i&0\\K_i&C_i&D_i\\0&L_i&0
\end{matrix}\right),
$$
where $K_i$ and $L_i$ are compact operators for $i=1, \ldots, n$. 

By Proposition \ref{p4} there exist quasi-nilpotent operators $Q$ and $R$
and compact operators $H_1, \ldots, H_n, M_1, \ldots, M_n$ such that $K_i=H_iQ$ and
$L_i =RM_i$  for $i=1,\ldots,n$.

Define operators $Q'_1, Q'_2$ and $\{S'_i\}_{i=1}^n$ on $\cH$  by 
$$Q'_1=\left(\begin{matrix}
    0&I&0\\ 0&0&I\\ R&0&0
\end{matrix}\right),\quad
Q'_2=\left(\begin{matrix}
    0&I&0\\ 0&0&I\\ Q&0&0
  \end{matrix}\right)\text{ and }
S'_i=\left(\begin{matrix}
    M_i&0&0\\A_i&0&0\\ C_i&D_i&H_i
  \end{matrix}\right) .$$
Thus $VT_iV^{-1}=Q'_1S'_iQ'_2$ for all $i=1, \ldots, n$. Since 
${Q'}_1^3= {\rm diag} \, (R, R, R)$, the operator $Q'_1$ is quasi-nilpotent. 
The same holds for the operator $Q'_2$. 

Let $Q_1=V^{-1}Q'_1V$, $Q_2=V^{-1}Q'_2V$ and
$S_i=V^{-1}S'_iV\quad(i=1.\dots,n)$. 
Clearly, the operators $Q_1$ and $Q_2$ are quasi-nilpotent,
and $T_i=Q_1S_iQ_2$ for all $i=1,\dots,n$.

Since there exists an orthonormal sequence
$\{h_k\}$ in $\cH$ such that $\|S'_ih_k\|\to 0\quad (k\to\infty)$
for all $i=1,\dots,n$, we have 
$\Bigl(\sum_{i=1}^n S^*_iS_i\Bigr)V^{-1}h_k\to 0$ and so $0\in
\sigma_e(\sum_{i=1}^{n} S^*_iS_i)$. 
Similarly it is possible to show that 
$0 \in \sigma_e(\sum_{i=1}^{n} S_iS^*_i)$.

This completes the proof of (c).
\end{proof}

\section{Case of nilpotent operators}

In the last section we prove the nilpotent analogs of Theorems \ref{two}
and \ref{c10}. We first characterize operators that are products of two nilpotent
operators.

\begin{theorem}\label{n1}
An operator $T \in \cBH$ is a product of two nilpotent operators (with index of
nilpotency
at most $3$) if and only if $\dim \ker T=\dim \ker T^* = \infty$.
\end{theorem}

\begin{proof}
Suppose that we can write $T=MN$, where $M$ and $N$ are nilpotent operators. Since
$\ker N \subset \ker T$ and $\ker M^* \subset \ker T^*$, and 
$\dim \ker N = \dim \ker M^* = \infty$, we have $\dim \ker T=\dim \ker T^* = \infty$.

Conversely, suppose that $\dim \ker T=\dim \ker T^* = \infty$. 
One can verify that we can choose a decomposition of $\cH$ as a direct sum of 
infinite-dimensional subspaces $\cH_1$, $\cH_2$ and $\cH_3$ such that 
$\cH_1 \subseteq \ker T$ and $\cH_3  \subseteq \ker T^*$. 
Then the matrix of $T$ relative to this decomposition is of the form
$$\left(\begin{matrix}
0&A&B \\ 0&C&D\\ 0&0&0
\end{matrix}\right).$$
Using the same argument as in the proof of Proposition \ref{p1} we can also achieve
that $B=0$.
Define operators $M$ and $N$ on $\cH$ by
$$M=
\left(\begin{matrix}
0&0&A\\I&0&C\\0&0&0
\end{matrix}\right)
\quad \text{ and } \quad
N=
\left(\begin{matrix}
0&0&D\\0&0&0\\0&I&0
\end{matrix}\right).
$$
It is easy to verify that $T=MN$ and that $M$ and $N$ are nilpotent operators with
index of nilpotency at most $3$.
\end{proof}

We now consider the nilpotent analog of Theorem \ref{c10}.
Given operators $T_1, \ldots, T_n$ on $\cH$, we are seeking nilpotent operators 
$N_1$ and $N_2$ and operators $S_1, \ldots, S_n$ on $\cH$ such that 
$T_i=N_1 S_i N_2$ for all $i=1, \ldots,n$. 
It is easy to see that the necessary condition for these factorizations is
infinite-dimensionality 
of both $\ker (\sum_{i=1}^n T_iT^*_i)$ and $\ker (\sum_{i=1}^n T^*_iT_i)$.
Note that $\ker (\sum_{i=1}^n T_iT^*_i)=\bigcap_{i=1}^n\ker T_i^*$
and $\ker (\sum_{i=1}^n T^*_iT_i)=\bigcap_{i=1}^n\ker T_i$.
The following result shows that this condition is also sufficient. 

\begin{theorem}\label{n2}
Let $T_1, T_2, \ldots, T_n$ be operators in $\cBH$ such that
$$ \dim \ker \left( \sum_{i=1}^n T_iT^*_i \right) = 
   \dim \ker \left( \sum_{i=1}^n T^*_iT_i \right) = \infty . $$
Then the following assertions hold:
\begin{enumerate}
\item[(a)] 
There exist nilpotent operators $N,N_1, N_2, \ldots, N_n$ on $\cH$ 
with index of nilpotency at most $3$ such that $T_i=N N_i N$ for all $i=1,2, \ldots,
n$.
\item[(b)] 
There exist nilpotent operators $N_1$ and  $N_2$ on $\cH$ with index of nilpotency 
at most $3$ and operators $S_1, \ldots, S_n$ on $\cH$ 
such that $T_i=N_1 S_i N_2$ for all $i=1,2, \ldots, n$, and the
products $N_1 S_i$ and $S_i N_2$ are nilpotent operators for all  $i=1,2, \ldots, n$.
\end{enumerate}
\end{theorem}

\begin{proof}
We can choose a decomposition of $\cH$ as a direct sum of infinite-dimensional
subspaces 
$\cH_1$, $\cH_2$ and $\cH_3$ such that  
$\cH_1 \subseteq \ker (\sum_{i=1}^n T^*_iT_i)$ and 
$\cH_3 \subseteq \ker (\sum_{i=1}^n T_i T_i^*)$. 
Since $\ker (\sum_{i=1}^n T^*_iT_i) \subseteq \ker T^*_jT_j = \ker T_j$ and 
$\ker (\sum_{i=1}^n T_iT^*_i) \subseteq \ker T_jT^*_j =\ker T^*_j$ for all $j=1,2,
\ldots, n$,
the operators $\{T_i\}_{i=1}^n$ are of the form
$$ T_i = \left(\begin{matrix}
0&A_i&B_i \\ 0&C_i&D_i\\ 0&0&0
\end{matrix}\right).
$$
As in the proof of Theorem \ref{c10} we can also achieve that $B_i = 0$ for all 
$i=1, \ldots, n$. 

In order to show part (a), we define operators $N$ and $\{N_i\}_{i=1}^n$ on $\cH$ by
$$ N=\left(\begin{matrix}
                0&I&0\\ 0&0&I\\ 0&0&0
           \end{matrix}\right) \quad
   \text{ and }
   N_i=\left(\begin{matrix}
                    0&0&0\\A_i&0&0\\ C_i&D_i&0
             \end{matrix}\right) . $$
Clearly, all of them are nilpotent operators with index of nilpotency at most $3$. 
Since $T_i=N N_i N$ for all $i=1,2, \ldots, n$, the proof of (a) is complete.

For part (b) we define nilpotent operators $N_1$ and  $N_2$ on $\cH$ 
(with index of nilpotency at most $3$) and operators $\{S_i\}_{i=1}^n$ on $\cH$ by
$$ 
N_1= \left(\begin{matrix}
                0&0&I\\I&0&0\\0&0&0
           \end{matrix}\right), \  
N_2= \left(\begin{matrix}
                0&0&I\\0&0&0\\0&I&0
              \end{matrix}\right)
\quad \text{ and } \quad
S_i=
\left(\begin{matrix}
D_i&0&C_i\\ 0&0&0\\ 0&0&A_i
\end{matrix}\right) . 
$$
It is easy to check that $T_i=N_1 S_i N_2$ for all $i=1,2, \ldots, n$, and 
products $N_1 S_i$ and $S_i N_2$ are nilpotent operators for all  $i=1,2, \ldots, n$.
This completes the proof of the theorem.
\end{proof}

\vspace{5mm}

{\bf\sc Acknowledgement.} The research was supported by  
grants No. 04-2003-04 and No. BI-CZ/03-04-4, Programme Kontakt of Czech and Slovenian
Ministries of Education. 
The first and the third author were supported in part 
by the Ministry of Higher Education, Science and Technology.
The second author was also partially
supported by grant No. 201/03/0041 of GA \v CR and by the
Institutional Research Plan AV0Z 10190503.

\bibliographystyle{amsplain}

\end{document}